\documentclass{amsart}
\usepackage{amsmath,amssymb,amsthm,bm}
\usepackage{tikz,caption,subcaption}
\usepackage{genyoungtabtikz}
\usepackage{ytableau,tikz,varwidth}
\usetikzlibrary{calc}
\usepackage{here}
\usepackage{cases}
\usepackage{dsfont}
\def\bbbz{\mathbb Z}

\def\cch{\mathcal H}
\def\ccsc{\mathcal {SC}}
\def\ccdd{\mathcal {DD}}
\def\ccp{\mathcal P}

\def\det{{\rm det}}

\def\sgn{{\rm sign}}

\DeclareMathOperator{\I}{I} 

\renewcommand{\sp}{\mathrm{sp}}
\DeclareMathOperator{\oo}{so}

\definecolor{RoyalPurple}{rgb}{0.47, 0.32, 0.66}

\numberwithin{equation}{section}
\newtheorem{thm}{Theorem}
\newtheorem{lem}{Lemma}
\newtheorem{proposition}[thm]{Proposition}
\newtheorem{definition}[thm]{Definition}

\theoremstyle{remark}

\newtheorem{ex}{Example}

\usepackage{mathtools, stmaryrd}
\usepackage{xparse} 
\DeclarePairedDelimiterX{\Iintv}[1]{\llbracket}{\rrbracket}{\iintvargs{#1}}
\NewDocumentCommand{\iintvargs}{>{\SplitArgument{1}{,}}m}
{\iintvargsaux#1} %
\NewDocumentCommand{\iintvargsaux}{mm} {#1\mkern1.5mu..\mkern1.5mu#2}

\usepackage[margin=1.45in]{geometry}

\begin{document}

\title[Combinatorial interpretations of the Macdonald identities for affine root systems]{Combinatorial interpretations of the Macdonald identities for affine root systems}


\author{David Wahiche}\address{Univ. Lyon, Universit\'e Claude Bernard Lyon 1, UMR 5208, Institut Camille Jordan, France}
\email{wahiche@math.univ-lyon1.fr}



\maketitle


\begin{abstract}
We explore some connections between vectors of integers and integer partitions seen as bi-infinite words. This methodology enables us to give a combinatorial interpretation of the Macdonald identities for affine root systems of the seven infinite families in terms of symplectic and special orthogonal Schur functions. From these results, we are able to derive $q$-Nekrasov--Okounkov formulas associated to each family. Nevertheless we only give results for types $\tilde{C}$ and $\tilde{C}^{\vee}$, and give a sketch of the proof for type $\tilde{C}$.
\end{abstract}


\section{Introduction and notations}

Formulas involving hook length abound in combinatorics and representation theory. Between 2006 and 2008, using various methods coming from representation theory \cite{We}, gauge theory \cite{NO} and combinatorics \cite{Ha}, several authors proved the so-called Nekrasov--Okounkov formula which can be stated as follows:
\begin{equation}\label{NOdebut}
\sum_{\lambda\in\ccp}T^{\lvert \lambda\rvert}\prod_{h\in\cch(\lambda)}\left(1-\frac{z}{h^2}\right)=\prod_{k\geq 1}\left(1-T^k\right)^{z-1}.
\end{equation}
Here $T$ is a formal variable, $z \in\mathbb{C}$, $\ccp$ is the set of integer partitions and $\cch(\lambda)$ is the multiset of hook lengths of the partition $\lambda$.

This formula does not only cover the generating series for $\ccp$ obtained by setting $z=0$ in \eqref{NOdebut}: it actually gives a connection between powers of the Dedekind $\eta$ function and integer partitions.
Among generalizations of \eqref{NOdebut} that can be found in the literature, a $(q,t)$-extension was proved by Rains--Warnaar \cite{RW}, by using refined skew Cauchy-type identities for Macdonald polynomials. This result was also obtained independently by Carlsson--Rodriguez-Villegas \cite{CRV} by means of vertex operators and the plethystic exponential. As mentioned in \cite{RW}, the special case $q=t$ is a reformulation of a result due to Dehaye--Han \cite{HD} and Iqbal--Nazir--Raza--Salem \cite{INRS} which reads as follows:
\begin{equation}\label{Hande}
\sum_{\lambda\in\ccp}T^{\lvert \lambda\rvert}\prod_{h\in\cch(\lambda)}\frac{(1-uq^h)(1-u^{-1}q^h)}{(1-q^h)^2}=\prod_{k,r\geq 1}\frac{(1-uq^rT^k)^r(1-u^{-1}q^rT^k)^r}{(1-q^{r-1}T^k)^r(1-q^{r+1}T^k)^r}.
\end{equation}
\noindent Here $T$ and $q$ are formal variables and $u \in \mathbb{C}$. Note that taking $u=q^z$ and letting $q\rightarrow 1$ in \eqref{Hande} yields \eqref{NOdebut}, although it is not immediate for the product side. 

Methods used by Han to prove \eqref{NOdebut} and Dehaye--Han for \eqref{Hande} both start from a specialization of the Macdonald formula for affine root systems specialized in type $\tilde{A}$. However, one needs an unspecialized Macdonald identity for type $\tilde{A}$ to get \eqref{Hande}. This can be found for instance in \cite{RS}, where Rosengren and Schlosser give a proof of Macdonald identities for the seven infinite affine root systems with elliptic determinantal evaluation (see also Stanton's reformulation in \cite{Stanton}). The next step in \cite{HD} to prove \eqref{Hande} uses new combinatorial notions such as exploded tableaux and a $V_t$-coding adapted from Garvan--Kim--Stanton \cite{GKS}. These techniques are close to the methodology presented here. However their extensions to other types seem complicated.

\begin{figure}[t]
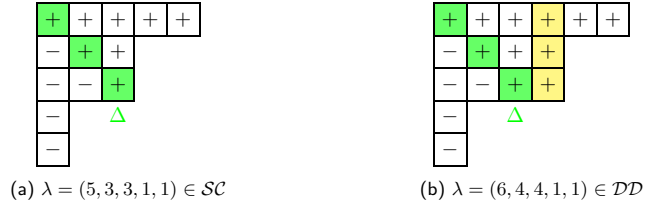

\centering
\scalebox{0.8}{
\begin{subfigure}[t]{.45\textwidth}
\centering

\begin{ytableau}
 *(green!60)+  &+ & + &+ & +\\
 - &*(green!60)+ &+ & \none \\
 - &-  &*(green!60)+&\none \\
  -&\none &\none[\textcolor{green}{\Delta}] &\none  \\
 - &\none &\none &\none \\
\end{ytableau}
\caption{$\lambda=(5,3,3,1,1)\in\ccsc$}

\label{fig:varsc}
\end{subfigure}
\begin{subfigure}[t]{.5\textwidth}
\centering
\begin{ytableau}
 *(green!60)+  &+ & + &*(yellow!60)+&+ & +\\
 - &*(green!60)+ &+ &*(yellow!60)+ & \none \\
 - &-  &*(green!60)+&*(yellow!60)+ &\none \\
  -&\none &\none[\textcolor{green}{\Delta}] &\none & \none \\
 - &\none &\none &\none & \none \\
\end{ytableau}
\caption{$\lambda=(6,4,4,1,1)\in\ccdd$}
\label{fig:vardd}
\end{subfigure}}
\caption{Self-conjugate and doubled distinct partitions filled with $\varepsilon$}
\label{fig:var}
\end{figure}

A \textit{partition} $\lambda$ of a positive integer $n$ is a nonincreasing sequence of positive integers $\lambda=(\lambda_1,\lambda_2,\dots,\lambda_\ell)$ such that $\lvert \lambda \rvert := \lambda_1+ \dots+\lambda_\ell = n$. The $\lambda_i$'s are the \textit{parts} of $\lambda$, the number $\ell$ of parts being the \textit{length} of $\lambda$, denoted by $\ell(\lambda)$. Each partition can be represented by its Ferrers diagram, which consists in a finite collection of boxes arranged in left-justified rows, with the row lengths in non-increasing order. The \textit{Durfee square} of $\lambda$ is the maximal square fitting in the Ferrers diagram. Its diagonal $\Delta$ will be called the main diagonal of $\lambda$. It is of size $d=d_\lambda:=\max(s | \lambda_s\geq s)$. Let us introduce a signed statistic $\varepsilon$ already appearing in \cite{MP}. For a box $s$ of $\lambda$ of coordinates $(i,j)$, \textit{$\varepsilon_{s}$} is defined as $-1$ if $s$ is strictly below the main diagonal of the Ferrers diagram and as $1$ otherwise, as depicted in Figure \ref{fig:var}. The partition $\lambda'=(\lambda_1',\lambda_2',\dots,\lambda_{\lambda_1}')$ is the \textit{conjugate} of $\lambda$, where $\lambda_j'$ denotes the number of boxes in the column $j$. 

Recall that the \textit{hook length} of $s$, denoted by $h_s$, is the number of boxes $v$ such that either $s=v$, or $v$ lies strictly below (respectively to the right) of $s$ in the same column (respectively row). For any $t\in \mathbb{N}^{*}$, the multiset of all hook lengths that are congruent to $0 \pmod t$ is denoted by $\mathcal{H}_t(\lambda)$. Note that $\mathcal{H}(\lambda)=\mathcal{H}_1(\lambda)$. A partition $\omega$ is a \textit{$t$-core} if $\cch_t(\omega)=\varnothing$. For any $\mathcal{A}\subset\ccp$, let $\mathcal{A}_{(t)}$ be the subset of elements of $\mathcal{A}$ that are $t$-cores. For example, the only $2$-cores are the ``staircase'' partitions $(k,k-1,\dots,1)$ for any $k\in \mathbb{N}^*$.

An integer partition $\lambda$ is \textit{self-conjugate} if its Ferrers diagram is symmetric along the main diagonal. Let $\mathcal{SC}$ be the set of self-conjugate partitions. The set of doubled distinct partitions, denoted by $\ccdd$, is that of all partitions $\lambda$ with Durfee square of size $d_\lambda$ such that $\lambda_i=\lambda_i'+1$ for all $i\in\lbrace 1,\dots,d\rbrace$. In Figure \ref{fig:varsc} for instance, $\lambda=(5,3,3,1,1)\in\ccsc$ has its main diagonal $\Delta$ shaded in green while the main diagonal of $\lambda=(6,4,4,1,1)\in\ccdd$ is shaded in green in Figure \ref{fig:vardd}. The strip shaded in yellow corresponds to the boxes added to a self-conjugate partition to obtain a $\ccdd$ partition. These subsets of partitions arise when one expresses the Weyl denominator formula for types $B,C$ and $D$ \cite[p.$79$]{Macbook} and have been of particular interest within the work  of P\'etr\'eolle \cite{MP} where two Nekrasov--Okounkov type formulas for $\ccsc$ and $\ccdd$ are derived. For instance, he proves the following $\ccdd$ Nekrasov--Okounkov type formula, coming from the basic specialization of the Macdonald identity specialized for type $\tilde{C}$, which stands that for a formal variable $T$ and any complex number $z$:
\begin{equation}\label{formulemathias}
\sum_{\lambda\in\ccdd} (-1)^{d_\lambda}T^{\lvert\lambda\rvert}\prod_{\substack{s\in\lambda\\h_s\in\cch(\lambda)}}\left(1-\frac{2z+2}{h_s\varepsilon_{s}}\right)=\prod_{k\geq 1}\left(1-T^k\right)^{2z^2+z}.
\end{equation}

%


Let $a$ and $T$ be formal variables. Recall that the $T$-Pochhammer symbol is defined as
\begin{equation*}
\displaystyle (a;T)_0=1, \quad (a;T)_\infty = \prod_{j\geq 0} (1-aT^j),\quad\text{and}\quad  \displaystyle (a_1,\dots,a_n)_\infty =(a_1;T)_\infty \dots (a_n;T)_\infty.
\end{equation*}

We denote by $\textsf{sp}$ (respectively $\textsf{so}$) the classical symplectic (respectively odd orthogonal) Schur function (see for instance \cite{FH}). Let $t$ be a strictly positive integer, set \linebreak ${\bf x}:=(x_1,\dots,x_t)$ and let us introduce:
$$ K_T(t,{\bf x})=\prod_{1\leq i<j\leq t}\left(Tx_i x_j,Tx_i^{-1} x_j^{-1},Tx_ix_j^{-1},Tx_jx_i^{-1};T\right)_\infty.$$

The goal of this extended abstract is to investigate combinatorially the connections between all the Macdonald identities and the Nekrasov--Okounkov formulas. According to Corollary 6.2 in \cite{RS} (see also \cite{Mac,Stanton}), the Macdonald identity for type $\tilde{C}_t$ is as follows:
\begin{multline}
\Delta_C(\mathbf{x})\left(T;T\right)_\infty^t \prod_{i=1}^t \left(Tx_i^2,Tx_i^{-2};T\right)_\infty K_T(t,\mathbf{x})=\sum_{{\bf m}\in\bbbz^t}\sum_{\sigma\in S_t}\sgn(\sigma)
\prod_{i=1}^t x_i^{(2t+2)m_i}\\\times T^{2(t+1)\binom{m_i}{2}+(t+1)m_i}
\left((x_iT^{m_i})^{\sigma(i)-t-1}-(x_iT^{m_i})^{t+1-\sigma(i)}\right),\label{eq:macdostart}
\end{multline}
where $\Delta_C(\mathbf{x})=\prod_{1\le i \le t}x_i^{-t}(1-x_i^2)\prod_{1\le i<j\le t}(x_j-x_i)(1-x_ix_j)$.
 We will focus on types $\tilde{C}_t$ (denoted $\tilde{C}_t^{(1)}$ in \cite{Kac}), $\tilde{C}_t^{\vee}$ (denoted $\tilde{D}_{t+1}^{(2)}$ in \cite{Kac}) here, but the method is the same for the other types. By computing the Littlewood decomposition to partitions seen as a bi-infinite sequences of ``$0$'' and ``$1$'' (see Section \ref{lit} for precise definitions and properties), the quadratic form which is to the exponent of $T$ in \eqref{eq:macdostart} can be interpreted as half the weight of a doubled distinct $(2t+2)$-core partition using the ideas of Garvan--Kim--Stanton \cite{GKS}. Introducing the notion of $V_{g,t}$-coding (see Definition \ref{def:vcoding}) which can be thought of as the last indices $\pmod{g}$ of letters ``$0$" in the bi-infinite sequences, and $H_+:=\lbrace h_s<g,\varepsilon_s=1\rbrace$, one can then reinterpret the right-hand side of\eqref{eq:macdostart} as follows:
\begin{thm}\label{prop:MacdotypeC}
Set $t\in\mathbb{N}^*$. The Macdonald identity for type $\tilde{C}_t$ can be rewritten as follows:
\begin{equation}\label{eq:macdoc}
\sum_{\omega\in\ccdd_{(2t+2)}} (-1)^{d_{\omega}+\lvert H_+\rvert}T^{\lvert\omega\rvert/2}\sp_{\mu}({\bf x})=
\left(T\right)_\infty^t K_T(t,\mathbf{x})\prod_{i=1}^t \left(Tx_i^2,Tx_i^{-2};T\right)_\infty ,
\end{equation}

where $\mathbf{v}$ is the $V_{g,t}-$coding corresponding to $\omega$ (see Definition \ref{def:vcoding}) and $\mu\in\ccp$ such that
$\mu_i:=v_i+i-2t-2$ for all $1\leq i \leq t$.
\end{thm}

Similarly we get the following:
\begin{thm}
\label{prop:Macdoothertypes}
Set $t\in\mathbb{N}^*$. The Macdonald identity for type $\tilde{C}_t^{\vee}$ can be rewritten as follows:
\begin{multline*}
\sum_{\omega\in\ccsc_{(2t)}} (-1)^{\lvert H_+\rvert+\lvert H_+\cap \Delta\rvert+d_\omega}T^{\lvert\omega\rvert/2}\oo_{\mu}({\bf x})\\
=\left(T^{1/2};T^{1/2}\right)_\infty\left(T;T\right)_\infty^{t-1}
K_T(t,{\bf x})\prod_{i=1}^t \left(T^{1/2}x_i,T^{1/2}x_i^{-1};T^{1/2}\right)_\infty
\end{multline*}
where $\mathbf{v}$ is the $V_{g,t}-$coding corresponding to $\omega$ (see Definition \ref{def:vcoding}) and $\mu\in\ccp$ is such that
$\mu_i:=v_i+i-2t$ for all $1\leq i \leq t$.
\end{thm}



As a consequence of our results, we can prove the following $q$-analogues of Nekrasov--Okounkov type identities.
\begin{thm}\label{NOC}
For formal variables $T$, $q$ and any complex number $u$, we have:
\begin{multline}\label{eq:noc}
\sum_{\lambda\in\ccdd} (-u)^{d_\lambda}T^{\lvert \lambda\rvert/2}\prod_{s\in\lambda}\frac{1-u^{-2\varepsilon_s}q^{h_s}}{1-q^{h_s}}\prod_{s\in\Delta}\frac{1+uq^{h_s/2}}{1+u^{-1}q^{h_s/2}}\\=
\prod_{m,r\geq 1}\frac{1+uq^{r-1}T^m}{1+u^{-1}q^{r}T^m}\frac{\left(1-u^{-2}q^{r+2}T^m\right)^{r-\lfloor r/2\rfloor}\left(1-u^2q^{r-1}T^m\right)^{r-\lfloor r/2\rfloor}}{\left(1-q^{r}T^m\right)^{r-\lfloor r/2\rfloor}\left(1-q^{r+1}T^m\right)^{r-\lfloor r/2\rfloor}}
\end{multline}
\begin{multline}\label{eq:nosc}
\sum_{\lambda\in\ccsc} (-1)^{d_\lambda}T^{\lvert \lambda\rvert/2}\prod_{s\in\lambda}\frac{1-u^{-2\varepsilon_s}q^{2h_s}}{1-q^{2h_s}}\prod_{s\in\Delta}\frac{1-uq^{h_s}}{1-u^{-1}q^{h_s}}\\=
\prod_{m,r\geq 1}\frac{\left(1-T^{m/2}\right)}{\left(1-T^{m}\right)}\frac{1-u^{-1}q^{2r-1}T^{m/2}}{1-uq^{2r-1}T^{m/2}}\frac{\left(1-u^{-2}q^{2(r+1)}T^m\right)^{r-\lfloor r/2\rfloor}\left(1-u^2q^{2r}T^m\right)^{r-\lfloor r/2\rfloor}}{\left(1-q^{2r}T^m\right)^{r+1-\lfloor (r+1)/2\rfloor}\left(1-q^{2(r+1)}T^m\right)^{r-\lfloor r/2\rfloor}}.
\end{multline}
%

\end{thm}

Note that taking $u=q^z$ and letting $q\rightarrow 1$ in \eqref{eq:noc} gives \eqref{formulemathias}, while in \eqref{eq:nosc} it yields a new Nekrasov--Okounkov type formula. Actually all the  Macdonald identities in the Appendix of \cite{Mac} can be derived from specializations and limits of Theorems \ref{prop:MacdotypeC} and \ref{prop:Macdoothertypes} and their analogues for other types.
In Section \ref{lit}, we introduce a way of computing the Littlewood decomposition with words, giving an explicit connection between $t$-cores and vectors of integers. In Section \ref{hooks}, we use this description of the Littlewood decomposition to characterize the product of hook lengths of subsets of $t$-cores such as $\ccdd_{(2t+2)}$. We derive Theorem \ref{thm:DDpair} but we only give the lemmas necessary to its proof. In the last section, we sketch the proofs of Theorems \ref{prop:MacdotypeC} and \ref{NOC} for type $\tilde{C}$ only.
\section{Combinatorial properties of the Littlewood decomposition on certain subsets of partitions} \label{lit}

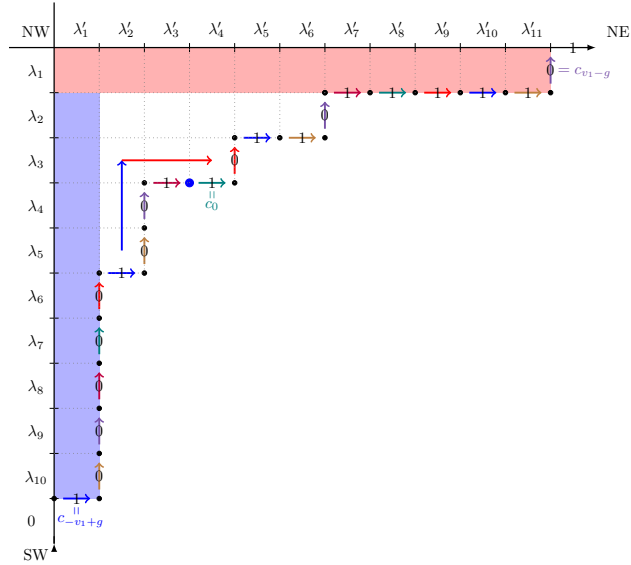
\begin{figure}[t]
\centering
\scalebox{0.6}{
\begin{tikzpicture}
    [
        dot/.style={circle,draw=black, fill,inner sep=1pt},
    ]
\fill [red!30] (0,0) rectangle (11,-1);
\fill [blue!30] (0,-1) rectangle (1,-10);
\foreach \x in {0,...,1}{
    \node[dot] at (\x,-10){ };
}

\draw[->,very thick,blue!100] (0.2,-10) -- (0.2+.6,-10);
\draw[->,very thick,blue!100] (1.2,-5) -- (1.2+.6,-5);
\draw[->,very thick,blue!100] (4.2,-2) -- (4.2+.6,-2);
\draw[->,very thick,blue!100] (9.2,-1) -- (9.2+.6,-1);

\draw[->,very thick,brown] (1,-9.8) -- (1,-9.8+.6);
\draw[->,very thick,brown] (2,-4.8) -- (2,-4.8+.6);
\draw[->,very thick,brown] (5.2,-2) -- (5.2+.6,-2);
\draw[->,very thick,brown] (10.2,-1) -- (10.2+.6,-1);

\draw[->,very thick,RoyalPurple] (1,-8.8) -- (1,-8.8+.6);
\draw[->,very thick,RoyalPurple] (2,-3.8) -- (2,-3.8+.6);
\draw[->,very thick,RoyalPurple] (6,-1.8) -- (6,-1.8+.6);
\draw[->,very thick,RoyalPurple] (11,-0.8) -- (11,-0.8+.6);

\draw[->,very thick,purple] (1,-7.8) -- (1,-7.8+.6);
\draw[->,very thick,purple] (2.2,-3) -- (2.2+.6,-3);
\draw[->,very thick,purple] (6.2,-1) -- (6.2+.6,-1);

\draw[->,very thick,teal] (1,-6.8) -- (1,-6.8+.6);
\foreach \x in {3.2}
    \draw[->,very thick,teal] (\x,-3) -- (\x+.6,-3);
\draw[->,very thick,teal] (7.2,-1) -- (7.2+.6,-1);

\draw[->,very thick,red] (1,-5.8) -- (1,-5.8+.6);
\draw[->,very thick,red] (4,-2.8) -- (4,-2.8+.6);
\draw[->,very thick,red] (8.2,-1) -- (8.2+.6,-1);

\foreach \y in {5,...,9}
	\node[dot] at (1,-\y){};
\node[dot] at (4,-3){};
\foreach \y in {3,...,5}
    \node[dot] at (2,-\y){};

\foreach \x in {4,...,6}
    \node[dot] at (\x,-2){};
\foreach \x in {6,...,11}
    \node[dot] at (\x,-1){};

\foreach \x in {1,...,11}
    \draw (\x,-.1) -- node[above,xshift=-0.4cm,yshift=1mm] {$\lambda_{\x}'$} (\x,+.1);

\node[above,xshift=0.5cm,yshift=1mm] at (12,0) {NE};
\node[above,xshift=-4mm,yshift=1mm] at (0,0) {NW};

\foreach \y in {1,...,10}
    \draw (.1,-\y) -- node[above,xshift=-4mm,yshift=0.2cm] {$\lambda_{\y}$} (-.1,-\y);
\node[below,xshift=-4mm] at (0,-11) {SW};
\node at (-0.5,-11+0.5){0};
\node at (11+0.5,0){1};

\foreach \y in {6,...,10}
	\node at (1,-\y+0.5){0};

\foreach \y in {4,...,5}
    \node at (2,-\y+0.5){0};

\node at (4,-3+0.5){0};
\node at (6,-2+0.5){0};
\node at (11,-0.5){0};

\node at (0.5,-10) {1};
\node at (1.5,-5) {1};

\foreach \x in {2,...,3}
    \node at (\x+0.5,-3){1};
\foreach \x in {4,...,5}
    \node at (\x+0.5,-2){1};
\foreach \x in {6,...,10}
    \node at (\x+0.5,-1){1};
\node at (0.5,-10.25){\textcolor{blue}{$\shortparallel$}};
\node at (0.5,-10.5){\textcolor{blue}{\phantom{o}$c_{-v_1+g}$}};

\node at (3.5,-3.25){\textcolor{teal}{$\shortparallel$}};
\node at (3.5,-3.5){\textcolor{teal}{$c_{0}$}};

\node at (11.75,-0.55){\textbf{\textcolor{RoyalPurple}{$=c_{v_1-g}$}}};

\node[dot] at (6,-1){};  
\draw[->,thick,-latex] (0,-11) -- (0,-11);
\draw[thick] (0,-11) -- (0,1);
\draw[->,thick,-latex] (-1,0) -- (12,0);
\node[circle,draw=blue,fill=blue,inner sep=0pt,minimum size=5pt] at (3,-3){};

\foreach \y in {6,...,9}
	\draw[dotted,gray] (0,-\y)--(1,-\y);
\draw[dotted,gray] (0,-5)--(2,-5);
\draw[dotted,gray] (0,-4)--(2,-4);
\draw[dotted,gray] (0,-3)--(3,-3);
\draw[dotted,gray] (0,-2)--(5,-2);
\draw[dotted,gray] (0,-1)--(6,-1);

\draw[dotted,gray] (1,-6)--(1,0);
\draw[dotted,gray] (2,-5)--(2,0);
\draw[dotted,gray] (3,-3)--(3,0);
\draw[dotted,gray] (4,-3)--(4,0);
\draw[dotted,gray] (5,-2)--(5,0);
\draw[dotted,gray] (6,-1)--(6,0);
\foreach \x in {7,...,10}
	\draw[dotted,gray] (\x,-1)--(\x,0);



\draw[->,very thick, blue](1.5,-4.5)--(1.5,-2.5);
\draw[->,very thick, red](1.5,-2.5)--(3.5,-2.5);

\end{tikzpicture}}
\caption{$\omega=(11,6,4,2,2,1,1,1,1,1)\in\ccdd_{(6)}$ and its binary correspondence}
\label{fig:wordddcore}
\end{figure}

In this section, we use the formalism of Han--Ji \cite{HJ}. Let $\partial \lambda$ be the border of the Ferrers diagram of $\lambda$. 
Encode the walk along the border from the South-West to the North-East as depicted in Figure \ref{fig:wordddcore}: take ``$0$'' for a vertical step and ``$1$'' for a horizontal step. This yields a $0/1$ sequence denoted \textbf{$s(\lambda)$}. This resulting word over the $\lbrace 0,1\rbrace$ alphabet has infinitely many ``$0$'''s at the beginning (respectively ``$1$'''s at the end), is indexed by $\bbbz$, and written $(c_i)_{i\in\mathbb{Z}}$.  
 
This writing as a sequence is not unique (since for any $k$ sequences, $(c_{k+i})_{i\in\mathbb{Z}}$ define the same partition), hence the necessity to set the index $0$ uniquely (to ensure this encoding is bijective). To tackle that issue, we set the index $0$ when the number of ``$0$'''s to the right of that index is equal to the number of ``$1$'''s to the left. In other words, the number of horizontal steps along $\partial\lambda$ corresponding to a ``$1$'' of negative index in $(c_i)_{i\in\bbbz}$ must be equal to the number of vertical steps corresponding to ``$0$'''s of nonnegative index in $(c_i)_{i\in\bbbz}$ along $\partial\lambda$. The delimitation between the letter of index $-1$ and that of index $0$ is called the \textit{median} of the word, marked by a $\mid$ symbol. The size of the Durfee square is then equal to the number of ``$1$'''s of negative index. Hence the application $s$ bijectively associates a partition to the word:
\begin{align*}
s(\lambda)=(c_i)_{i\in\mathbb{Z}}=\left(\ldots c_{-2}c_{-1}|c_0c_1c_2\ldots\right), \intertext{where $c_i\in\lbrace 0,1\rbrace$ for any $i\in\bbbz$, and such that}
\#\{i\leq-1,c_i=1\}= \#\{i\geq0,c_i=0\}<\infty.
\end{align*}

\begin{lem}\label{lem:indices}
This application maps bijectively a box $s$ of hook length $h_s$ of the Ferrers diagram of $\lambda$ to a pair of indices $(i_s, j_s)\in\bbbz^2$ of $s(\lambda)$ such that $i_s<j_s$, $c_{i_s}=1$ and $c_{j_s}=0$, $j_s-i_s=h_s$.
\end{lem}  
Lemma \ref{lem:durf} below allows to characterize the position of a box.
\begin{lem}(\cite[Lemma 2.1]{Wmult})\label{lem:durf}
Set $\lambda\in\ccp$ and $s(\lambda)$ its corresponding word. Let $s$ be a box of the Ferrers diagram of $\lambda$. Let $(i_s,j_s)\in\bbbz^2$ be the indices in $s(\lambda)$ associated with $s$. Then $s$ is a box strictly above the main diagonal in the Ferrers diagram of $\lambda$ if and only if $|i_s|\leq |j_s|$.
\end{lem}
The map below is often called the Littlewood decomposition (see \cite{GKS,HJ} for instance).

\begin{definition}\label{defphi}{\em
Let $t \geq 2$ be an integer and consider:\\
$$\begin{array}{l|rcl}
\Phi_t: & \mathcal{P} & \to & \mathcal{P}_{(t)} \times \mathcal{P}^t \\
& \lambda & \mapsto & (\omega,\nu^{(0)},\ldots,\nu^{(t-1)}),
\end{array}$$
where if $s(\lambda)=\left(c_i\right)_{i\in\bbbz}$, then for all $k\in\lbrace 0,\dots,t-1\rbrace$ one has $\nu^{(k)}:=s^{-1}\left(\left(c_{ti+k}\right)_{i\in\bbbz}\right)$. The tuple $\underline{\nu}=\left(\nu^{(0)},\ldots,\nu^{(t-1)}\right)$ is the $t$-quotient of $\lambda$, denoted by \textit{$quot_t(\lambda)$}, while $\omega$ is the $t$-core of $\lambda$, denoted by \textit{$core_t(\lambda)$}.}
\end{definition}
Obtaining the $t$-quotient is straightforward from $s(\lambda)=\left(c_i\right)_{i\in\bbbz}$: we just look at subwords with indices congruent to the same values modulo $t$. The sequence $10$ within these subwords are replaced iteratively by $01$ until the subwords are all the infinite sequence of ``$0$'''s before the infinite sequence of ``$1$'''s (in fact it consists in removing all rim hooks in $\lambda$ of length congruent to $0\pmod t$). Then $\omega$ is the partition corresponding to the word which has the subwords $\pmod t$ obtained after the removal of the $10$ sequences.

For example, if we take $\lambda = (4,4,3,2) \text{ and } t=3$, then $s(\lambda)=\ldots \color{red}{0} \color{blue}{0} \color{green}{1} \color{red}{1} \color{blue}{0}\color{green}{1} \color{black}| \color{red}{0} \color{blue}{1} \color{green}{0} \color{red}{0} \color{blue}{1} \color{green}{1}\color{black}\ldots$
\begin{align*}
\begin{array}{rc|rcl}
s\left(\nu^{(0)}\right)=\ldots \color{red} 001 \color{black}| \color{red}001\color{black}\ldots& &s\left(w_{0}\right)=\ldots \color{red} 000 \color{black}| \color{red}011\color{black}\ldots,\\
 s\left(\nu^{(1)}\right)=\ldots \color{blue} 000 \color{black}| \color{blue}111\color{black}\ldots& \longmapsto& s\left(w_{1}\right)=\ldots \color{blue} 000 \color{black}| \color{blue}111\color{black}\ldots , \\
 s\left(\nu^{(2)}\right)=\ldots \color{green} 011 \color{black}| \color{green}011\color{black}\ldots & & s\left(w_{2}\right)=\ldots \color{green} 001 \color{black}| \color{green}111\color{black}\ldots .
\end{array}\\
\end{align*}
Thus $s(\omega)=\ldots \color{red}{0} \color{blue}{0} \color{green}{0} \color{red}{0} \color{blue}{0}\color{green}{1} \color{black}| \color{red}{0} \color{blue}{1} \color{green}{1} \color{red}{1} \color{blue}{1} \color{green}{1}\color{black}\ldots $, and
$$
quot_3(\lambda)=\left(\nu^{(0)},\nu^{(1)},\nu^{(2)}\right)=\left((1,1),\varnothing,(2)\right),\ core_3(\lambda)= \omega=(1)
$$

One might therefore see $t$-core partitions as partitions whose $t$-quotient in the Littlewood decomposition is empty. Let $\omega$ be a partition, it is a $t$-core if and only if $quot_t(\omega)=\left(\varnothing,\dots,\varnothing\right)$. This is equivalent to say that all subwords $\pmod t$ are of the form $\ldots 0011\ldots$, which is, an infinite sequence of ``$0$'''s and then an infinite sequence of ``$1$'''s. For any $i\in\lbrace 0,\dots,t-1\rbrace$ let us define $n_i:=\min\lbrace k \in\bbbz\mid c_{i+kt}=1\rbrace$. Each $n_i$ corresponds to the index of the first ``$1$'' in the subword of $s(\omega)$ whose index is congruent to $i \pmod t$. Recall that the word $s(\omega)$ has as many ``$1$'''s of negative index as ``$0$'''s of positive index. This is equivalent to require $\sum_{i=0}^{t-1}n_i=0$, hence there is a natural bijective map between $\omega\in \ccp_{(t)}$ and $\left(n_i\right)_{i\in \lbrace 0,\dots,t-1\rbrace}\in\bbbz^t$ such that $\sum_{i=0}^{t-1}n_i=0$.

For example, if we take $\omega = (4,2) \text{ and } t=3$, then 
\begin{align*}
\begin{array}{rc|lc}
& &s\left(w_{0}\right)=\ldots \color{red} 000 &\color{black}|\color{red}\underbrace{001}_{n_0=2}1\color{black}\ldots,\\
 s\left(\omega\right)=\ldots \color{red}{0} \color{blue}{0} \color{green}{0} \color{red}{0} \color{blue}{1}\color{green}{1} \color{black}| \color{red}{0} \color{blue}{1} \color{green}{1} \color{red}{0} \color{blue}{1} \color{green}{1}\color{red}{1}\color{blue}{1} \color{green}{1}\color{black} \ldots& \longmapsto& s\left(w_{1}\right)=\ldots \color{blue} 00\underbrace{1}_{n_1=-1}&\color{black}| \color{blue}111\color{black}\ldots , \\
  & & s\left(w_{2}\right)=\ldots \color{green} 00\underbrace{1}_{n_2=-1} &\color{black}| \color{green}111\color{black}\ldots .
\end{array}\\
\end{align*}

In \cite{Johnson}, Johnson uses the fermionic viewpoint of partitions (which is the same as the one described above) to prove that this bijection is actually the one used by Garvan--Kim--Stanton in \cite{GKS}. We reformulate here what is written in \cite[Section 2]{Johnson} in terms of index of words: let $\lambda$ be a partition and $t$ be a strictly positive integer. Abaci correspond exactly to $t$-subwords of $s(\lambda)$ with fixed residue $\pmod{t}$ while the $n_i$'s as defined in Bijection $2$ in \cite{GKS} correspond to the charge of the $i$-th runner on the abaci. This implies in particular that if $\omega$ is a $t$-core partition with corresponding word $s(\omega)=\left(c_i\right)_{i\in\bbbz}$ and if we set $\phi(\omega):=\left(n_0,\dots,n_{t-1}\right)$, then $n_i=\min\lbrace k\in\bbbz \mid c_{kt+i}=1\rbrace$. Moreover, we have:
\begin{equation}
|\omega|=\frac{t}{2}\sum_{i=0}^{t-1} n_i^2+\sum_{i=0}^{t-1} i n_i.
\label{eq:gks}
\end{equation}


%

\section{Hook length product of $t$-core partitions}\label{hooks}
The aim of this section is to introduce the material required to prove Theorem \ref{thm:DDpair} which is an enumerative result on hook length products. The latter is crucial to get Theorem \ref{NOC}, as explained in the introduction (see also Section \ref{sketch} for more details). This extended abstract only focuses on results for $\tilde{C}$ and $\tilde{C}^{\vee}$ that are stated independently in order to avoid characteristic functions (that complicate their statement). However note that analogous statements exist for all seven infinite families. As for technical results allowing to prove Theorem \ref{thm:DDpair}, they are stated only in the $\ccdd_{(2t+2)}$-case, once again to avoid characteristic functions and technicality they bring. The cases $\ccsc_{(2t)}$, $\ccsc_{(2t+1)}$ and $\ccdd_{(2t+1)}$ are obtained on the same way, uniform formulations and proofs can be found in \cite{Wmac}. Theorem \ref{thm:DDpair} involves $V_{g,t}$-codings: these are a particular kind of vector of integers (see Definition \ref{def:vcoding}) that can be associated to any $\omega \in \ccdd_{(2t+2)}$ (see Proposition \ref{prop:core_vcoding}).

%
%
%

Let $t\in \mathbb{N}^*$. The end of Section \ref{lit} associates a vector of integers $\phi(\omega)$ to any $t$-core $\omega$. 
Following \cite{CRV}, this vector can be ordered ``naturally'' (meaning that this order is inherited by that of the residue $\pmod{t}$) by setting:
\begin{equation}\label{eqn:phi}
n_k:=\left\lfloor \frac{\lambda_i-i}{t}\right\rfloor+1, \quad i=\min\lbrace \nu \mid \lambda_\nu \equiv k \pmod t\rbrace.
\end{equation}

The Littlewood decomposition, when restricted to $\ccdd$, also has interesting properties and can be stated as follows (\cite{GKS,MP} for instance):
\begin{numcases}
{\lambda\in\ccdd \mapsto} \left(\omega,\underline{\tilde{\nu}}\right)\in\ccdd_{(t)}\times\ccdd\times\ccp^{(t-1)/2}\quad\text{if $t$ is odd,}\notag  \\
\left(\omega,\underline{\tilde{\nu}},\mu\right)\in\ccdd_{(t)}\times\ccdd\times\ccp^{(t-2)/2}\times\ccsc\quad\text{if $t$ is odd.}\label{todd}
\end{numcases}

This symmetrical behaviour of $\ccdd_{(2t+2)}$ partitions yields some additional conditions on the associated vector of integers. These have already been studied by Garvan--Kim--Stanton \cite{GKS} but are stated here in a slightly different way. Let $\omega\in \ccdd_{(2t+2)}$ and $\phi(\omega)=(n_0,n_1,\dots,n_{2t+1})$ be as above. Symmetries of $\ccdd$ ensure that $n_0=0$, and $n_i=-n_{2t+2-i}$ for all $i \in \lbrace 1,\dots,2t+1\rbrace$, which in particular implies that $n_t=0$. By \eqref{eq:gks}, we then have:
\begin{align}
\lvert \omega \rvert &=(t+1)\sum_{i=0}^{2t+1}n_i^2+\sum_{i=0}^{2t+1}in_i=2\left((t+1)\sum_{i=0}^{t-1}n_i^2+\sum_{i=0}^{t-1}(i-t-1)n_i)\right)\label{eq:gkddpair}
\end{align}
\noindent So with this formalism one can recover the vector of integers given in \cite{GKS} from the Littlewood decomposition together with the word decomposition. When developing \eqref{eq:macdostart} with respect of powers of $T$, half of the previous quadratic form appear to the exponent.
\begin{ex}
Figure \ref{fig:wordddcore} above illustrates this for a $\ccdd$ when $t=6$. The arrows are sorted in six different colors, each of them corresponding to a fixed residue $\pmod 6$ of the index of the corresponding word of $\omega$.
The word corresponding to $\omega$ writes as follows:
\begin{center}
$s(\omega)=\cdots \color{blue}{0} \color{brown}{0}  \color{RoyalPurple}{0} \color{purple}{0} \color{teal}{0} \color{purple}{0}\color{blue}{1} \color{brown}{0}  \color{RoyalPurple}{0} \color{purple}{0} \color{teal}{0} \color{purple}{0} \color{blue}{1} \color{brown}{0}  \color{RoyalPurple}{0} \color{purple}{1} \color{black}|\color{teal}{1} \color{purple}{0} \color{blue}{1} \color{brown}{1}  \color{RoyalPurple}{0} \color{purple}{1} \color{teal}{1} \color{purple}{1} \color{blue}{1} \color{brown}{1}  \color{RoyalPurple}{0} \color{purple}{1} \color{teal}{1} \color{purple}{1}
\color{blue}{1} \color{brown}{1}  \color{RoyalPurple}{1} \color{purple}{1} \color{teal}{1} \color{purple}{1}\color{black}\cdots.$
\end{center}

By extracting the subwords of fixed residue $\pmod 6$, we obtain:
\begin{align*}
&s\left(w_{0}\right)=\cdots \color{teal} 000 \color{black}| \color{teal}111\color{black}\cdots, \quad
s\left(w_{1}\right)=\cdots \color{purple} 000 \color{black}| \color{purple}011\color{black}\cdots \quad
s\left(w_{2}\right)=\cdots \color{blue} 011 \color{black}| \color{blue}111\color{black}\cdots\\
&s\left(w_{3}\right)=\cdots \color{brown} 000 \color{black}| \color{brown}111\color{black}\cdots, \quad
s\left(w_{4}\right)=\cdots \color{RoyalPurple} 000 \color{black}| \color{RoyalPurple}001\color{black}\cdots, \quad
s\left(w_{5}\right)=\cdots \color{purple} 001 \color{black}| \color{purple}111\color{black}\cdots
\end{align*}

so that $(\color{teal}{n_0},\color{purple}{n_1},\color{blue}{n_2},\color{brown}{n_3},\color{RoyalPurple}{n_4},\color{purple}{n_5}\textcolor{black}{)=(}\color{teal}{0},\color{purple}{1},\color{blue}{-2},\color{brown}{0},\color{RoyalPurple}{2},\color{purple}{-1} \textcolor{black}{)\in\bbbz^6.}$

The properties of symmetry of any $\ccdd$ partition along its main diagonal yields the following restriction, as mentioned in \cite{GKS}:
\begin{center}
$ \omega \in\ccdd_{(6)}\longleftrightarrow (\color{purple}{1},\color{blue}{-2}\textcolor{black}{)\in\bbbz^2.}$
\end{center}
\end{ex}
Similarly, the symmetrical behaviour of $\ccsc_{(2t)}$ partitions yields additional conditions on the associated vector of integers: let $\omega\in \ccsc_{(2t)}$ and $\phi(\omega)=(n_0,n_1,\dots,n_{2t-1})$, then $n_i=-n_{2t-1-i}$ for all $i \in \lbrace 0,\dots,2t-1\rbrace$, in particular $n_t=0$. By \eqref{eq:gks}, we have:
\begin{align}
\lvert \omega \rvert &=(2t)/2\sum_{i=0}^{2t-1}n_i^2+\sum_{i=0}^{2t-1}in_i=\sum_{i=0}^{t-1}((2t)n_i^2+(2(i-t)-1)n_i).\label{eq:gkscimpair}
\end{align}

The ordered vector of integers $\phi(\omega)$ still lacks of some properties to be the combinatorial tool one needs to prove Theorems \ref{prop:MacdotypeC}, \ref{NOC} and \ref{thm:DDpair}. These properties are the one satisfied by vectors of integers called $V_{g,t}$-codings (defined in \cite{Wmac}). The goal here is to attach a $V_{g,t}$-coding to any $\phi(\omega)$. Before going any further, let us motivate Definition \ref{def:vcoding} below, with the following observation: making use of properties of $q$-series together with the Weyl group action on the affine root systems, Stanton provides a new proof of the Macdonald identities \cite{Stanton}. Sums that appear in these identities involve vectors of integers, so that the existence of $\phi$, which is a bijection between partitions and vectors of integers, leads us to rephrase Stanton's work. From this perspective, \cite[Proposition 3.7]{Stanton} that focuses on type $\tilde{C}_t$ establishes why only partitions in $\ccdd_{(2t+2)}$ are to be considered on type $\tilde{C}_t$: they are actually the only ones whose associated vector of integers contribute in a non-trivial way to the sum in the Macdonald identity for the type $\tilde{C}_t$.  

As this extended abstract only focuses on results for types $\tilde{C}_t$ and $\tilde{C}_t^{\vee}$, Definition \ref{def:vcoding} below restricts to these types. Both the general definition (of \cite{Wmac}) and Definition \ref{def:vcoding} involve a parameter $g$ of the affine root system defined by the equality $M=g\Lambda$ in \cite[p.$134$]{Mac}. In particular $g = 2t+2$ if the affine algebra considered is of type $\tilde{C}_t$ and $ g = 2t$ if it is of type $\tilde{C}_t^{\vee}$. 

\begin{definition}\label{def:vcoding}
Let $t \in \mathbb{N}^*$ and $g \in \mathbb{N}^*$ such that $t\leq g$, then set $\lambda \in \ccp$ and let $s(\lambda)=(c_k)_{k\in\bbbz}$ be its corresponding binary word. For $i\in\lbrace 0,\dots,g-1\rbrace$, define $\beta_i:=\max\lbrace (k+1)g+i\mid c_{kg+i}=0\rbrace$. Let $\sigma:\lbrace 1,\dots,g\rbrace\rightarrow \lbrace 0,\dots,g-1\rbrace$ be the unique permutation such that $\beta_{\sigma(1)}>\dots>\beta_{\sigma(g)}$. The vector $\mathbf{v}:=(\beta_{\sigma(1)},\dots,\beta_{\sigma(t)})$ is called the $V_{g,t}$-coding corresponding to $\lambda$.
\end{definition}


Note that a notion of $V_{t}$-coding introduced by Dehaye--Han in \cite{HD} to connect Macdonald identities for type $\tilde{A}_{t-1}$ and $t$-core partitions corresponds to the $V_{t,t}$-coding. We can now state Theorem \ref{thm:DDpair} and provide the key technical results required to prove it.
 
 \begin{thm}\label{thm:DDpair}
Set $t$ a positive integer and $g=2t+2$. Let $\omega\in \ccdd_{(2t+2)}$ and $\mathbf{v}\in\bbbz^{t}$ its associated $V_{g,t}$-coding, and set $r_i=v_i-t-1$ for any $i\in\lbrace1,\dots,t\rbrace$. Then we have
\begin{equation*}\label{eq:ddpoids}
\lvert \omega\rvert = \frac{1}{g}\sum_{i=1}^{t} r_i^2-\frac{(g/2-1)(g-1)}{12},
\end{equation*}
and setting $\alpha_i(\omega):=\#\lbrace u\in\omega, h_u=g-i, \varepsilon_u=1\rbrace$, and for any function $\tau:\bbbz\rightarrow F^{\times}$, where $F$ is a field, we also have
\begin{equation*}\label{eq:thmdd}
\prod_{s\in\omega}\frac{\tau(h_s-\varepsilon_s g)}{\tau(h_s)}=\prod_{i=1}^{g-1}\left(\frac{\tau(-i)}{\tau(i)}\right)^{\alpha_i(\omega)}
\prod_{i=1}^{t}\frac{\tau(r_i)}{\tau(i)}\prod_{1\leq i<j\leq t} \frac{\tau(r_i-r_j)}{\tau(j-i)} \frac{\tau(r_i+r_j)}{\tau(g-i-j)}.
\end{equation*}

%
\end{thm}

To show Theorem \ref{thm:DDpair} one needs to associate bijectively a vector of integers to any partition. For technical and conceptual reasons (respectively facilitating inductive proofs and interpreting the determinants as Schur functions) the involved vectors need to be ordered. The $V_{g,t}$-codings introduced so far are natural candidate, but one still need to show that they correspond bijectively to partitions. This is the purpose of Proposition \ref{prop:core_vcoding}.
\begin{proposition}\label{prop:core_vcoding}
Let $t$ be a positive integer. Any $\lambda$ in one of the sets $\ccdd_{(2t+2)}$ and $\ccsc_{(2t)}$ is in bijective correspondence with its $V_{g,t}$-coding, where $g$ is the index of the corresponding set.
\end{proposition}
 
Let $s(\omega)$ be the word corresponding to $\omega$ via the word decomposition. The $V_{g,t}$-coding associated to $\omega$ by Proposition \ref{prop:core_vcoding} is immediately given when reading $s(\omega)$ from right to left: for instance $v_1$ is the index of the first ``$1$'' in the subword $\pmod{g}$ of $s(\omega)$ that contains the last ``$0$''.
Lemma \ref{lem:maxhook} allows to prove Theorem \ref{thm:DDpair} by induction on the length of the Durfee square of the partition. Following the same philosophy as in Lemmas \ref{lem:indices} and \ref{lem:durf}, it enumerates the boxes in the largest hook by means of their pair of indices. Ultimately, and as illustrated in Figure \ref{fig:wordddcore} in the red and blue shaded areas, let $s$ be a box of $\lambda$ in the first hook. Then $\varepsilon_s=1$, respectively $\varepsilon_s=-1$, implies $j_s=v_1-g$, respectively $i_s=-v_1+g$ and $j_s<v_1-g$. 

We define the \textit{$g$-intervals} as as $\I_{m,M}^{g,+}:=\lbrace k\in \bbbz \mid m\leq k<M, k\equiv m \pmod g\rbrace$ and $\I_{m,M}^{g,-}:=\lbrace l\in \bbbz \mid m<l\leq M, l\equiv M \pmod g\rbrace$. This notion is of particular interest in our case ever since we are trying to enumerate hook lengths with a fixed residue $\pmod g$. The key ingredient of the proof of Theorem \ref{thm:DDpair} is that $g$-intervals allow us to transform products involving generic functions $\tau$ into telescopic factors.

\begin{lem}\label{lem:maxhook}
Set $\omega\in \ccdd_{(2t+2)}$ so that $g = 2t+2$ and let $(v_1,\dots,v_{t})\in\bbbz^t$ be its associated $V_{g,t}$-coding. Then the biggest hook of $\omega$, denoted by $H_1$, corresponds to the collection of boxes $\cch_{1,+}\cup\cch_{1,-}$ where $\cch_{1,+}$ is the set of indices of boxes $s$ in the first hook such that $\varepsilon_s=1$:
\begin{equation*}
\cch_{1,+}=\I_{-v_1+g,v_1-g}^{g,+}  \cup \I_{0,v_1-g}^{g,+} \cup\I_{g/2,v_1-g}^{g,+}
 \displaystyle\bigcup_{i=2}^{t}\left(\I_{v_i,v_1-g}^{g,+}\cup \I_{-v_i+g,v_1-g}^{g,+}\right),
 \end{equation*}
and $\cch_{1,-}$ is the set of indices of boxes $s$ in the first hook such that $\varepsilon_s=-1$:
\begin{equation*}
\cch_{1,-}= \I_{-v_1+g,v_1-2g}^{g,-} \cup\I_{-v_1+g,-g}^{g,-}\cup\I_{-v_1+g, -g/2}^{g,-} \displaystyle\bigcup_{i=2}^{t}\left(\I_{-v_1+g,v_i-g}^{g,-}\cup \I_{v_1-g,-v_i}^{g,-}\right).
\end{equation*}

\end{lem}

Finally we will also need the following lemma to apply Theorem \ref{thm:DDpair} to Macdonald identities specializations.
\begin{lem}\label{lem:signschur}
Let $\omega\in\ccdd_{(2t+2)}$ and $s(\omega)=(c_k)_{k\in\bbbz}$. Let $u_i=(2t+2)n_i+i$ the index of the first letter ``$1$'' in the subword $(c_{i+(2t+2)k})_{k\in\bbbz}$, where $(n_i)=\phi^{-1}(\omega)$.
Let $(v_i)_{i\in\lbrace 1,\dots t\rbrace}$ be a $V_{g,t}$-coding. Then $\sigma: \lbrace 1,\dots,t\rbrace\rightarrow\lbrace 0,\dots 2t+1\rbrace$ of Proposition \ref{def:vcoding} is such that that $v_i=u_{\sigma(i)}$ for $i \in \{1, \dots t\}$. We have that:
\begin{equation*}
\lvert H_+\rvert=\#\lbrace s\in\omega, h_s<2t+2, \varepsilon_s=1\rbrace\equiv d_\omega+\sgn(\sigma)\pmod{2}.
\end{equation*}
\end{lem}

\section{Sketch of proof of Theorems \ref{prop:MacdotypeC} and \ref{NOC}}\label{sketch}
 We first expand the right-hand side of \eqref{eq:macdostart} and we extract the terms in $T$. The power of $T$ corresponds to \eqref{eq:gkddpair}. By changes of variables and Lemma \ref{lem:signschur}, we are able to take $\omega=\phi^{-1}({\bf m})$ and derive Theorem \ref{prop:MacdotypeC}.



We now derive Theorem \ref{NOC} for $\ccdd$. In order to do so, we start by proving the equality setting $u=q^t$ for any strictly positive integer $t$ on both sides. By technical manipulations on products, the right hand side of \eqref{eq:noc} is exactly that of Theorem \ref{prop:MacdotypeC} when $x_i=q^i$ which completes the first step of the proof.

 The left-hand side of \eqref{eq:noc} can be obtained by setting $\tau(x)=1-q^x$ in Theorem \ref{thm:DDpair} and  multiplying the resulting expression by the product of hook lengths on the main diagonal $\Delta$. We prove that:
 \begin{equation}\label{eq:schurinter}
\sp_{\mu}(q,q^2,\dots,q^t)=(-1)^{\lvert H+\rvert} q^{(t+1)d_\lambda}\prod_{s\in\lambda}\frac{1-q^{h_s-(2t+2)\varepsilon_s}}{1-q^{h_s}}\prod_{s\in\Delta}\frac{1+q^{t+1+h_s/2}}{1+q^{-t-1+h_s/2}}.
 \end{equation}
The proof of the equality above is done by induction on $r_1=\max_{1\leq i\leq t}(r_i)$ as defined in Theorem \ref{thm:DDpair}. First remark that the product over elements of $\Delta$ on the right-hand side of \eqref{eq:schurinter} is a telescopic product equal to $\prod_{i=1}^t (1+q^{r_i})/(1+q^i)$. On the left-hand side of \eqref{eq:schurinter}, the quotient of determinants $\det(q^{ir_j}-q^{-ir_j})/\det(q^{i(t+1-j)}-q^{-i(t+1-j)})$ can be computed by noting that it is equal to $\det({(q^{r_j})}^{i}-{(q^{r_j})}^{-i})/\det({(q^{t+1-j})}^i-{(q^{t+1-j})}^{-i})$. Finally we conclude by checking that both sides of the equation verify the same initial value and the same induction property.

The last step of the proof is to check that both sides of \eqref{eq:noc} are Laurent polynomials in the variable $u$. An argument of polynomiality that \eqref{eq:noc} holds for any $u$ then allows one to conclude.

\subsection*{Acknowledgements}
The author would like to thank the anonymous referees from the FPSAC committee for their very enlightening comments to improve significantly the quality and the readability of this extended abstract. The author would also like to thank Benjamin Dupont, Marion Jeannin, Philippe Nadeau and Nicolas Ressayre for their help and the fruitful conversations they had. Any critical remark must be exclusively addressed to the author of this paper.

\bibliographystyle{plain}
\bibliography{macdo_ident}

\begin{thebibliography}{10}

\bibitem{CRV}
Erik Carlsson and Fernando Rodriguez~Villegas.
\newblock Vertex operators and character varieties.
\newblock {\em Adv. Math.}, 330, 2018.

\bibitem{HD}
Paul-Olivier Dehaye and Guo-Niu Han.
\newblock A multiset hook length formula and some applications.
\newblock {\em Discrete Math.}, 311(23-24), 2011.

\bibitem{FH}
William Fulton and Joe Harris.
\newblock {\em Representation theory}, volume 129 of {\em Graduate Texts in
  Mathematics}.
\newblock Springer-Verlag, New York, 1991.
\newblock A first course, Readings in Mathematics.

\bibitem{GKS}
F.~Garvan, D.~Kim, and D.~Stanton.
\newblock Cranks and {$t$}-cores.
\newblock {\em Invent. Math.}, 101(1), 1990.

\bibitem{Ha}
G.-N. Han.
\newblock The {N}ekrasov-{O}kounkov hook length formula: refinement, elementary
  proof, extension and applications.
\newblock {\em Ann. Inst. Fourier (Grenoble)}, 60(1), 2010.

\bibitem{HJ}
G.-N. Han and K.~Q. Ji.
\newblock Combining hook length formulas and {BG}-ranks for partitions via the
  {L}ittlewood decomposition.
\newblock {\em Trans. Amer. Math. Soc.}, 363(2), 2011.

\bibitem{INRS}
A.~Iqbal, S.~Nazir, Z.~Raza, and Z.~Saleem.
\newblock Generalizations of {N}ekrasov-{O}kounkov identity.
\newblock {\em Ann. Comb.}, 16(4), 2012.

\bibitem{Johnson}
P.~Johnson.
\newblock Lattice points and simultaneous core partitions.
\newblock {\em Electron. J. Combin.}, 25(3):Paper No. 3.47, 2018.

\bibitem{Kac}
Victor~G. Kac.
\newblock {\em Infinite-dimensional {L}ie algebras}.
\newblock Cambridge University Press, third edition, 1990.

\bibitem{Mac}
I.~G. Macdonald.
\newblock Affine root systems and {D}edekind's {$\eta $}-function.
\newblock {\em Invent. Math.}, 15, 1972.

\bibitem{Macbook}
I.~G. Macdonald.
\newblock {\em Symmetric functions and {H}all polynomials}.
\newblock Oxford Mathematical Monographs. The Clarendon Press, Oxford
  University Press, New York, second edition, 1995.
\newblock With contributions by A. Zelevinsky, Oxford Science Publications.

\bibitem{NO}
N.~A. Nekrasov and A.~Okounkov.
\newblock Seiberg-{W}itten theory and random partitions.
\newblock In {\em The unity of mathematics}, volume 244 of {\em Progr. Math.}
  Birkh\"{a}user Boston, Boston, MA, 2006.

\bibitem{MP}
M.~P\'{e}tr\'{e}olle.
\newblock {\em {Quelques d\'{e}veloppements combinatoires autour des groupes de
  Coxeter et des partitions d'entiers}}.
\newblock Theses, {Universit{\'{e}} Claude Bernard - Lyon I}, November 2015.

\bibitem{RW}
Eric~M. Rains and S.~Ole Warnaar.
\newblock A {N}ekrasov-{O}kounkov formula for {M}acdonald polynomials.
\newblock {\em J. Algebraic Combin.}, 48(1), 2018.

\bibitem{RS}
Hjalmar Rosengren and Michael Schlosser.
\newblock Elliptic determinant evaluations and the {M}acdonald identities for
  affine root systems.
\newblock {\em Compos. Math.}, 142(4), 2006.

\bibitem{Stanton}
Dennis Stanton.
\newblock An elementary approach to the {M}acdonald identities.
\newblock In {\em {$q$}-series and partitions ({M}inneapolis, {MN}, 1988)},
  volume~18 of {\em IMA Vol. Math. Appl.} Springer, New York, 1989.

\bibitem{Wmult}
David Wahiche.
\newblock Multiplication theorems for self-conjugate partitions.
\newblock {\em Comb. Theory}, 2(2):Paper No. 13, 2022.

\bibitem{Wmac}
David Wahiche.
\newblock Combinatorial versions of macdonald identities for affine root
  systems and its applications, work in progress.

\bibitem{We}
B.~W. Westbury.
\newblock Universal characters from the {M}acdonald identities.
\newblock {\em Adv. Math.}, 202(1), 2006.

\end{thebibliography}

\end{document}